\documentclass[a4paper,10pt,twocolumn]{article}
\usepackage[english]{babel}
\usepackage[utf8]{inputenc}
\usepackage[T1]{fontenc}
\usepackage[ruled,vlined]{algorithm2e}

\DecMargin{2.5mm}
\usepackage[top=1.5cm, left=1.5cm, right=1.5cm, bottom=1.5cm]{geometry}

\renewenvironment{abstract}{\bf\small {\em\ Abstract---}}{}

\usepackage{amsfonts,amssymb,amsmath,amsthm}
\usepackage{mathabx}
\usepackage{subfigure}
\usepackage{graphicx}
\usepackage[footnotesize]{caption}
\usepackage[dvipsnames]{xcolor}
\usepackage[numbers]{natbib}
\usepackage[colorlinks, linkcolor = red, citecolor=blue]{hyperref}
\usepackage{amsfonts,amssymb,amsmath,amsthm}
\usepackage{mathrsfs}

\newtheorem{lemma}{Lemma}
\newtheorem{theorem}{Theorem}

\DeclareMathOperator*{\argmin}{arg\,min}
\DeclareMathOperator{\R}{\mathbb{R}}
\DeclareMathOperator{\C}{\mathbb{C}}
\DeclareMathOperator{\N}{\mathbb{N}}
\DeclareMathOperator*{\rank}{rank}

\usepackage{empheq}
\title{Denoising and Completion of Structured Low-Rank Matrices via Iteratively Reweighted Least Squares}

\author{Christian K\"ummerle$^1$, Claudio Mayrink Verdun$^{1}$.\\
  \footnotesize $^1$Technische Universit\"at M\"unchen, Munich, Germany.  \\
  } \date{\empty} % no need for a date

\begin{document}

\maketitle

\begin{abstract} 
We propose a new Iteratively Reweighted Least Squares (IRLS) algorithm for the problem of completing or denoising low-rank matrices that are \emph{structured}, e.g., that possess a Hankel, Toeplitz or block-Hankel/Toeplitz structure. The algorithm optimizes an objective based on a non-convex surrogate of the rank by solving a sequence of quadratic problems.

Our strategy combines computational efficiency, as it operates on a lower dimensional generator space of the structured matrices, with high statistical accuracy which can be observed in experiments on hard estimation and completion tasks.

Our experiments show that the proposed algorithm \texttt{StrucHMIRLS} exhibits an empirical recovery probability close to 1 from fewer samples than the state-of-the-art in a Hankel matrix completion task arising from the problem of spectral super-resolution of badly separated frequencies.

Furthermore, we explain how the proposed algorithm for structured low-rank recovery can be used as preprocessing step for improved robustness in frequency or line spectrum estimation problems.
\end{abstract}

\section{Introduction}
\label{sec:introduction}
In this paper, we consider two related problems: first, let $X\in \C^{d_1 \times d_2}$ be a Hankel matrix with $X = (X_{i,j}) = x_{i+j-1}$ with some $x \in \C^{n}$, $n= d_1 +d_2 -1$, and consider the task to approximate $X$ by a low-rank Hankel matrix $\widehat{Z}$ such that
\begin{equation} \label{eq:Hankel_approx}
\widehat{Z} =\argmin\limits_{\substack{\C^{d_1 \times d_2} \ni Z \text{  is Hankel} \newline \\ \rank(Z) \leq r}} \|Z- X\|_{F(w)}^2,
\end{equation}
where $\|\cdot\|_{F(w)}$ is a suitable weighted Frobenius norm \cite{condat_hirabayashi,markovsky} and $r  < \min(d_1,d_2)$. Secondly, let $X= (X_{i,j})= (x_{i+j-1})$ be a Hankel matrix, $\Phi : \C^{n} \rightarrow \C^{m}$ a subsampling operator, $\mathcal{H}: \C^{n} \rightarrow\C^{d_1 \times d_2}, x \mapsto \mathcal{H}(x)= (\mathcal{H}(x)_{i,j})= (x_{i+j-1})$ be the linear Hankel operator and let
%The problem of recovering a low-rank matrix from incomplete observations, such as a given subset of its entries, is now ubiquitous in several fields of science \cite{davenport_romberg}, \cite{chen_chi_bigdata}. The huge amount of data generated  in some applications can be viewed as lying approximately in a low-dimensional subspace where some low-rank structure can be exploited.% \cite{candes_recht}. 

%Despite the fact that the general problem of finding the lowest rank matrix satisfying some measurement constraints is NP-hard, several algorithms were proposed to tackle the problem and, under certain assumptions, recovery guarantees have been derived for them.
%
%Among the algorithms to solve this problem, non-convex approaches are often preferred in practice because of their higher empirical recovery rate and their more efficient implementation \cite{chen_chi_bigdata}. Nevertheless, most of those methods do not take into account any other implicit information beyond low-rank assumptions.
\begin{equation} \label{eq_Hankel_completion}
\text{Find } \widehat{Z}=\mathcal{H}(\widehat{z}) \text{ s.t. } \widehat{z} =\argmin\limits_{z \in \C^{n}, \Phi(z)=\Phi(x)} \rank(\mathcal{H}(z)),
\end{equation}
the task of completing the Hankel matrix $X$ from a subset $\Phi(x)$ of its entries \cite{cai_wang_wei,cai_wang_wei2,ChenChi13,Ongie16}.

Problems as \eqref{eq:Hankel_approx} and \eqref{eq_Hankel_completion} arise from engineering applications as parallel MRI \cite{jin_lee_ye}, system identification \cite{LiuV10}, direction of arrival \cite{yang_li_stoica_xie} and seismic data interpolation \cite{gao_sacchi_chen}, since low-rankness of corresponding Hankel or block-Hankel matrices emerges from sparsity of the signal in a \emph{continuous} transform domain, generalizing in some sense the concept of \emph{sparsity} in discrete domains.

We note that for the problem of frequency estimation from samples of a sum of exponentials, also known as super-resolution \cite{candes_fernandez-granda,TangRecht13}, line spectral estimation \cite{tang_bhaskar_recht} or spectral compressive sensing \cite{ChenChi13}, low-rank properties of structured matrices play a significant role as well \cite{roy_kaylath,yang_xie_stoica}.

%\CK{Since finding a polynomial algorithm to minimize \eqref{eq:Hankel_approx} or \eqref{eq_Hankel_completion} exactly cannot be expected in general \cite{__}, the quest for algorithms with good approximate solutions or exact solutions under particular assumptions arises naturally.}
Due to the interplay of linear (block)-Hankel and non-convex low-rank structures, finding efficient computational approaches for \eqref{eq:Hankel_approx} and \eqref{eq_Hankel_completion} is challenging, especially for large-scale and multi-dimensional problems \cite{ChenChi13,jin_lee_ye}, where the quadratic growth of dimensionality from $x$ to $\mathcal{H}(x)$ is an %computational
issue.

In this work, we develop a new Iteratively Reweighted Least Squares (IRLS) algorithm called \emph{Structured HM-IRLS} tailored to the structured low-rank matrix estimation problems \eqref{eq:Hankel_approx} and \eqref{eq_Hankel_completion}. It builds on previous work on algorithms for the recovery of unstructured low-rank matrices \cite{Fornasier11,KS18,Mohan10}. % for the structured low-rank matrix recovery problem that overcomes some limitations from previous works, being scalable and robust to noise, enabling us to obtain state-of-the-art in terms of recovery results and also shedding light on the line spectral estimation problem. \vspace{-.3cm}

\section{Our Approach}
\label{sec:first-section}
To derive an optimization-based approach for the problems \eqref{eq:Hankel_approx} and \eqref{eq_Hankel_completion}, the idea is to replace the non-convex and non-smooth rank by a (still) non-convex \emph{logdet-surrogate}

\begin{equation*}  %\label{eq:def_schattenp} 
\|X\|_{0,\epsilon}^0 :=  \sum_{i=1}^{d}\log(\sigma_i (X)^2+  \epsilon^2) % \log \det (X X^* + \epsilon I_{d_1})  = 
\end{equation*}
for some $\epsilon > 0$ with $d= \min(d_1,d_2)$, if $\sigma_i(X)$ denotes the $i$-th singular value of $X$. Using this notation, if $X = \mathcal{H}(x)$, our algorithm can be seen as a sequential minimization of local quadratic upper bounds on  

\begin{equation}\label{eq:modelp} 
\mathcal{J}_{\lambda}(z,\epsilon):= \lambda \|\mathcal{H}(z)\|_{0,\epsilon}^0 +  \|\Phi(z)- y\|_2^2
\end{equation}
where $\lambda$ is a regularization parameter chosen such that $\lambda \approx 0$ for the completion problem \eqref{eq_Hankel_completion} and $\lambda > 0$ for the denoising problem \eqref{eq:Hankel_approx}, and $y \in \R^m$ is the data vector. For the denoising problem, $\Phi$ corresponds to the identity with $m = n$.

We summerize the steps of \texttt{StrucHMIRLS} in Algorithm 1, denoting by $\mathcal{T}_R$ be the rank-$R$ truncation operator. 
\vspace*{-3mm}
\begin{algorithm}[h]
	\DontPrintSemicolon 
	\KwIn{Map $\Phi: \mathbb{R}^{n} \rightarrow \mathbb{R}^m$, data vector $y \in \mathbb{R}^m$, rank estimate $R$, regularization parameter\hspace{-0.5mm} $\lambda > 0$,\hspace{-0.1mm} $\alpha < 1$.}
%	\KwOut{Sequence $(X^{(n)})_{n=1}^{n_0} \subset M_{d_1 \times d_2}$.}
	Initialize $k=0, \epsilon_0\!=\!\sigma_1(\mathcal{H}(\Phi^*(y))), W^{(0)}\!\! = \!\epsilon_0^2 \mathbf{I}_{n}\! \in \C^{n \times n}$. \\
	\lRepeat{$\|z^{(k)}-z^{(k-1)}\|_2/\|z^{(k)}\|_2 < \operatorname{tol}$.}{ $k = k+1$,
		\vspace*{-2mm}
	\begin{align}
			z^{(k)} &=\argmin\limits_{z \in \mathbb{C}^{n}}\;\; \langle z,W^{(k)} z\rangle_{\ell_2} + \frac{1}{2 \lambda} \|\Phi(z)-y\|_2^2, \nonumber\\ %\label{xmin} \\
			\epsilon_k &=\min\left(\epsilon_{k-1},||z^{(k-1)}-z^{(k)}||_2+\alpha^{k^2}\right), \nonumber\\ %\label{epsmin}\\
			\!\!\!\!\!\!\!\!W^{(k+1)}\!&=\!2 \mathcal{H}_{\operatorname{vec}}^*\Big[\Big(\mathcal{T}_R\big(\mathcal{H}(z^{(k)})\mathcal{H}(z^{(k)})^*\big) \oplus \nonumber\\
			& \mathcal{T}_R\big(\mathcal{H}(z^{(k)})^*\mathcal{H}(z^{(k)})\big)\Big) + \epsilon_k^2  \mathbf{I}_{d_1 d_2}\Big]^{-1}\mathcal{H}_{\operatorname{vec}}, \nonumber %\label{Wmin} \\
		\end{align}}
	\caption{\!\texttt{Structured harmonic mean\! IRLS \newline (StrucHMIRLS)}} 
\end{algorithm} 

%Let $\mathcal{T}_R$ be the rank-$R$ truncation operator. Pseudo-code for \texttt{StrucHMIRLS} is given in Algorithm 1. 
%Let $\mathcal{T}_R$ be the rank-$R$ truncation operator \CK{and $\mathcal{H}_{\operatorname{vec}}$ be the operator that maps to the canonical vectorization of the Hankel matrix.} 
Furthermore, $\mathcal{H}_{\operatorname{vec}}$ denotes the operator that maps to the canonical vectorization of the Hankel matrix, and   %In the notation of Algorithm 1, we use
 $\oplus$ is the Kronecker sum $A\oplus B= \mathbf{I}_{d_2} \otimes A + B \otimes \mathbf{I}_{d_1}$ of matrices $A \in \C^{d_1 \times d_1}$, $B \in \C^{d_2 \times d_2}$, while $\otimes$ denotes the standard Kronecker product. %\cite{KS18}.

We can show that the sequence $(z^{(k)},\epsilon_k)_{k \in \N}$ does never increase the value of the \emph{logdet} objective $\mathcal{J}_{\lambda}(z,\epsilon)$ from \eqref{eq:modelp}.
\begin{theorem}
The sequence $\big(\mathcal{J}_{\lambda}(z^{(k)},\epsilon_k) = \lambda \|\mathcal{H}(z^{(k)})\|_{0,\epsilon_k}^0 + \|\Phi(z^{(k)})- y\|_2^2\big)_{k\in \N}$ is non-increasing.
\end{theorem}

For the design of the algorithm, it is crucial that the quadratic upper bounds on $\mathcal{J}_{\lambda}(z^{(k)},\epsilon_k)$ are well-chosen, since the optimization landscape of $\mathcal{J}_{\lambda}$ is in general extremely non-convex. Our particular choice of the weight matrix $W^{(k+1)}$ in Algorithm 1 achieves precisely this, as it uses the information in both the column and the row space of $\mathcal{H}(z^{(k)})$. Our reweighting rule can be considered as an extension of the \emph{harmonic mean} weight matrix rule of \cite{KS18}, which corresponds to the harmonic mean of the left- and right-sided weight matrices used in the IRLS strategies of \cite{Fornasier11,Mohan10} in the unstructured case. 

We note that while an IRLS algorithm for structured low-rank matrix recovery has been already proposed in \cite{Ongie16}, the authors of \cite{Ongie16} do not optimize a rank surrogate of the Hankel matrix $\mathcal{H}(z)$  itself, but of a half-circulant extension which is not expected to be very low-rank even if $\mathcal{H}(z)$ is.

\vspace{-3mm}
\section{Computational efficiency}
\vspace{-1mm}
If an upper estimate $R$ of the rank $r$ of the (unknown) ground truth matrix $X= \mathcal{H}(x)$ is known, the steps of \texttt{StrucHMIRLS} can be implemented in a efficient way: To define the weight matrices, it is sufficient to use randomized SVDs \cite{Halko11} of order $R$, which can be calculated with complexity $O(n R \log(n))$ due to fast multiplication of vectors by a Hankel matrix $\mathcal{H}(z)$. Furthermore, we obtain the following lemma about fast matrix-vector multiplication by the matrix $W^{(k)}$ of Algorithm 1.
\begin{lemma}
The multiplication of $W^{(k)} \in \C^{n \times n}$ with a vector $v \in \C^n$ can be computed in $O(nR^2+ n R \log n)$ operations.
\end{lemma}
This can be used to update the $z^{(k)}$ efficiently using a \emph{conjugate gradient method} for solving the linear systems corresponding to the weighted quadratic problems.
\vspace{-3mm}
\section{Connection to Super-resolution}
\label{sec:superesolution}
\vspace{-2mm}
In certain applications such as seismology, radar and wireless communications, the goal is to recover signals that are specified by parameters in a continuous domain as, e.g., a mixture of r complex sinusoids. In this case, one can model the signal of interest such that $x_0(t)=\sum_{i=1}^r \alpha_i e^{i 2 \pi f_i t}$ with unknown frequencies $\{f_1, \dots, f_r\} \subset [0,1]$ and amplitudes $(\alpha_i)_{i=1}^r \in \mathbb{C}$. If we are given a vector $x=(x(0),\ldots,x(n-1)) \in\C^n$ of (noiseless) discrete samples of $x_0$, it can be seen that the corresponding Hankel matrix $\mathcal{H}(x)$ has rank $r$ if $d_1$ is chosen such that $r< \min(d_1, n-d_1+1)$, if the $r$ frequencies $f_i$ do not coincide \cite{Hua92,jin_lee_ye}.

In case of missing samples $x_T \in \C^{|T|}$ corresponding to $T \subset \{0, \dots, n-1\}$, $|T| = m$, or in case of noise on the samples, this low-rank structure can be exploited via \texttt{StrucHMIRLS}, to be combined with an algorithm as ESPRIT \cite{roy_kaylath} or others \cite{Hua92,schmidt} to obtain the frequencies in a two-stage procedure.

%Due to the Vandermonde decomposition of a Hankel matrix \cite{yang_xie_stoica}, one is able to establish a connection between the sum of a few exponentials with the low-rank Hankel matrices formed by the samples and then this can be solved exactly as a structured matrix completion problem, without assuming any discretization of the frequencies. Several techniques were proposed to solve it in a off-the-grid manner \cite{Hua92}, \cite{cadzow}, \cite{roy_kaylath}, \cite{schmidt}

% but most of them are very sensitive to noise or to minimum separation between the frequencies \cite{tang_bhaskar_recht}. 

%Once one have all the samples in a denoised form, we could then recover the frequencies and amplitudes by using, for example, the ESPRIT \cite{roy_kaylath} algorithm, which performs well in the absense of noise. Therefore in our work we used the HM-IRLS combined with ESPRIT to obtain superior performance in the line spectral estimation problem.
 \vspace{-3mm}
\section{Numerical Results}
\label{sec:numerics}
In a first experiment, we consider the completion of Hankel matrices $\mathcal{H}(x)$ from $m$ sample coordinates $T$ of $x=(x(0),\ldots,x(n-1))$ that are drawn uniformly at random, $n=127$ and $x(t)=\sum_{i=1}^r\alpha_i e^{(2\pi i f_i)t}$, where the $f_i$ and $c_i$ are sampled independently such that $f_i \sim \mathcal{U}([0,1]),|\alpha_i| = 1+10^{c_i}$, $c_i \sim \mathcal{U}([0,1])$, as in \cite{cai_wang_wei}. In Figure $1$, the empirical recovery probabilities averaged over 50 simulations for each pair of $m$ and $r$ are documented in comparison with the  of algorithms \cite{tang_bhaskar_recht,cai_wang_wei,cai_wang_wei2,ChenChi13,fang_wang_shen_li_blum}. \texttt{StrucHMIRLS} exhibits the best performance, with successful recovery already when $m \approx 2 r$, despite the fact that some frequencies $f_i$ will be very close if $r$ is not too small, which compromises the performance of, e.g., atomic norm minimization \cite{tang_bhaskar_recht}.  

\begin{figure}[h] %\label{fig_1}
\includegraphics[width=0.51\textwidth]{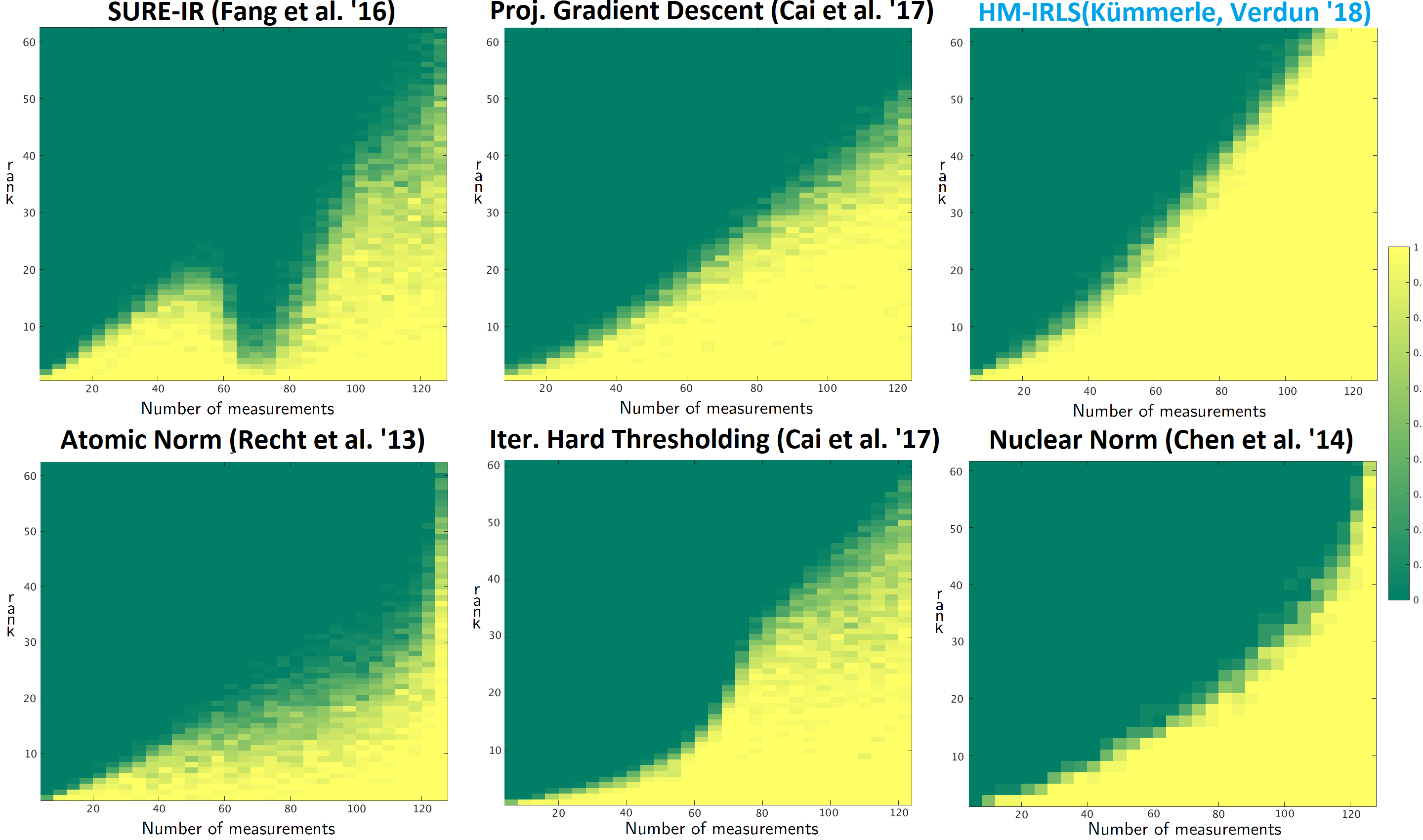} 
 \label{fig:1} \vspace{-5mm}\caption{Hankel matrix completion, $m$ measurements of vector $x \in \C^n$ with $n=123$. x-axis: number of measurements $m$, y-axis: model order $r$.}
\end{figure}
\vspace{-2mm}
Furthermore, we use the denoising variant of \texttt{StrucHMIRLS} in a second experiment to investigate its performance for frequency estimation under the presence of additive Gaussian noise on equispaced samples from a signal that is a sum of two frequencies located at $f_1= 0.35$ and $f_2=0.40$ (both with unitary amplitude), following Section VI of \cite{haro_vetterli}. 
\vspace{-2mm}
\begin{figure}[h]
\begin{center}
	\includegraphics[width=0.49\textwidth]{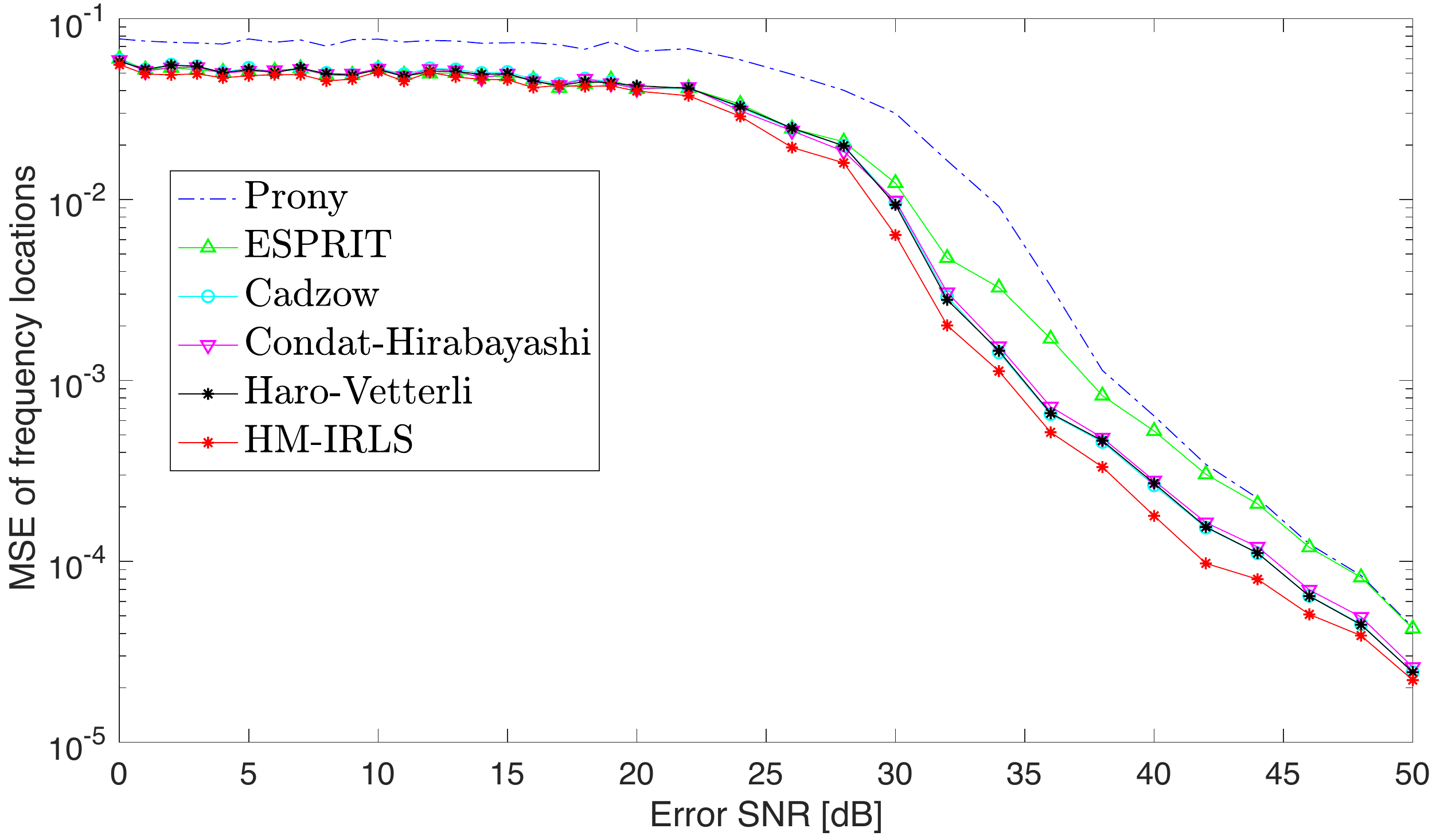}
\end{center}
\vspace{-5mm}
\caption{Frequency estimation experiment, $r=2$, $\alpha_1=\alpha_2=1$ and $f_1= 0.35$ and $f_2=0.40$.}
\end{figure}
\vspace{-2mm}

After denoising, we use ESPRIT to obtain the frequencies, which is arguably one of best algorithms for frequency estimation for low noise levels. As a comparison, we use the algorithms \cite{cadzow,condat_hirabayashi,haro_vetterli} (combined with ESPRIT for frequency retrieval, respectively), vanilla-ESPRIT \cite{roy_kaylath} and Prony's method. For our method, we choose the regularization parameter $\lambda$ according to an adaptive rule that uses the information of the model order $r=2$.

The results corresponding to an average over 500 independent noise realizations for each SNR value can be seen in Figure 2, and our method consistently obtains a lower MSE on the vector of frequencies $f=(f_1,f_2)$ than the competing methods across different noise SNRs. 

%the noise model and it can be shown that $\lambda=\sigma\sqrt{n\log n}$ is the optimal choice for a small asymptotic MSE \cite{tang_bhaskar_recht} in the case of atomic norm minimization. Since our go is to outperform this strategy, in this work it will also be adopt as a rule of thumb.

\bibstyle{abbrv}
\bibliography{Dissertation_Literature_instance2}

\end{document}